\documentclass[12pt]{amsart}
\usepackage{amssymb, amscd, amsmath}
\usepackage{latexsym}
\usepackage{epsf}
\usepackage{amssymb, epic,eepic,epsfig,amsbsy,amsmath,amscd}
\newtheorem{thm}{Theorem}[section]
\newtheorem{dfn}[thm]{Definition}
\newtheorem{Thm}{Theorem}

\newtheorem{lem}[thm]{Lemma}

\theoremstyle{definition}

\newtheorem{remark}{Remark}[section]
\def\1{{\rm1\mathchoice{\kern-0.25em}{\kern-0.25em}
        {\kern-0.2em}{\kern-0.2em}I}}

\newcommand{\one}{{ \rm \setlength{\unitlength}{1em}
\begin{picture}(0.75,1)
\put(0,0){$1$}\put(0.34,0){\line(0,1){0.65}}
\end{picture} }}

\newcommand{\lmn}[1]{\vadjust{\setbox1=\vtop{\hsize 25mm
\parindent=0pt\baselineskip=9pt
\rightskip=4mm plus 4mm#1}
\hbox{\kern-26mm\smash{\raise .5ex\box1}}}}
\input epsf.sty

\newcommand{\nc}{\newcommand}
\def\be#1\ee{\begin{equation}#1\end{equation}}
\nc{\bc}{\begin{center}} \nc{\ec}{\end{center}} \nc{\bb}{\mathbb}
\nc{\cal}{\mathcal} \nc{\frk}{\mathfrak} \nc{\N}{{\mathsf N}}
\nc{\K}{{\mathsf K}} \nc{\fk}{\mathbf{k}} \nc{\fn}{\mathbf{n}}
\nc{\fb}{\mathbf{b}}  \nc{\e}{\varepsilon} \nc{\ev}{{\rm{ev}}}

\def\ll{\left[\!\!\left[}
\def\rr{\right]\!\!\right]}

\theoremstyle{remark}

\def\Z{{\mathbb Z}}

\def\N{{\mathbb N}}
\def\C{{\mathcal C}}

\def\v8{\vskip 8pt}

\def\e{\varepsilon}

\def\la{\langle}
\def\ra{\rangle}

\def\R{\mathbb R}

\def\BiC{\operatorname{Ex\cal C}}

\def\Mor {{\operatorname{Mor}}}
\def\id {{\operatorname{id}}}

\def\deg{\operatorname{nes}}
\def\gr{\operatorname{deg}}

\def\Kh{\operatorname{Kh}}
\def\Cob{\operatorname{Cob}}
\def\Odd{\operatorname{OddCob}}
\def\NC{\operatorname{NesCob}}
\def\EC{\operatorname{EmbCob}}
\def\Nkho{\operatorname{NesKho}}
\def\mod{\operatorname{Mod}}

\def\ep{\epsilon}

\begin{document}

\title[Nested Khovanov Homology]
{On  link homology theories\\[0.2cm] from extended cobordisms}

\author{Anna Beliakova}
\address{Institut f\"ur Mathematik, Universit\"at Z\"urich,
 Winterthurerstrasse 190,
CH-8057 Z\"urich, Switzerland}
\email{anna@math.unizh.ch}

\author{Emmanuel Wagner}
\address{Institut de Math\'ematiques de Bourgogne, Universit\'e de Bourgogne, UMR 5584 du CNRS, BP47870, 21078 Dijon Cedex, France}

\email{emmanuel.wagner@u-bourgogne.fr}


\begin{abstract}
This paper is devoted to the study of algebraic structures leading to link
homology theories. The originally used structures of Frobenius algebra
and/or TQFT
are modified in two directions. First, we refine
 2--dimensional cobordisms by taking 
into account their embedding
into $\R^3$. Secondly, we extend the underlying cobordism category
to a 2--category, where the usual relations hold up to 2--isomorphisms.
The corresponding abelian 2--functor is called 
 an extended quantum field theory
(EQFT). We show that the  Khovanov homology, the nested Khovanov homology,
extracted by Stroppel and Webster from Seidel--Smith
construction,
and the odd Khovanov homology fit into this setting.
Moreover, we prove
that any EQFT based on a $\Z_2$--extension 
of the embedded cobordism category
which coincides with
Khovanov after reducing the coefficients modulo 2, gives rise
to a link invariant homology theory isomorphic to those of Khovanov.
\end{abstract}

\maketitle
\section*{Introduction}
In his influential paper \cite{Kho},  Khovanov constructed a link homology theory categorifying
the Jones polynomial. During few years,
this categorification was considered to be essentially unique, since the underlying
$(1+1)$ TQFT was known to be  determined by its Frobenius system and all
 rank two Frobenius systems were fully classified \cite{Kho1}.
However, in \cite{Odd} Ozsvath,  Rasmussen and Szabo
 came up with a new categorification of 
the Jones polynomial,  which agrees with  Khovanov's one
after reducing the coefficients modulo two. The underlying algebraic structure
of the odd Khovanov homology can not be described in terms of the
 Frobenius algebra.

This fact attracts attention again to the question of description and
classification of algebraic structures leading to link homology theories. 
In this paper, we provide an  evidence to the fact that the appropriate  algebraic structure 
is given by an extended quantum field theory (EQFT). 
A EQFT here is a 2--functor
from a certain (semistrict) monoidal 2--category of cobordisms, called 
an {\it extension},  to an abelian category. 
Given a cobordism category by specifying its generators and relations, the
 2--category is constructed by requiring the relations to be satisfied 
up to 2--isomorphisms.
Furthermore, such a 2--category is called
 an extension of the original cobordism category, if
 the automorphism group of any 1--morphism is trivial. 
A  simple example of an extension is a $\Z_2$--extension, where 
the 2--morphisms
are just plus or minus the identity.
 Notice that extensions can be defined for both strict and semistrict
monoidal 2--categories and the resulting EQFT will also be called 
strict and semistrict respectively.

The usage of the word ``extension'' in our setting 
is motivated by the fact that
after replacing  the original category by a group we will get a
usual extension of that group. Those extensions are
 classified by the second cohomology
classes of the group. Therefore, our approach can serve as a definition
for the second cohomology of a category.
A quite different notion of an 
extended topological field theory (ETFT) was introduced
and studied in \cite{SP}. 


In this paper, we construct extensions of the category of
2--dimensional cobordisms $\Cob$ and of the category  of embedded 2--cobordisms
modulo the unknotting relation $\NC$. In the first case, we 
recover the Khovanov and the odd Khovanov homologies, as strict (trivial) and 
semistrict extensions respectively. 
In the second case,
we construct so--called nested Khovanov homology, 
extracted by Stroppel and Webster \cite{CW} 
from the algebraic counterpart of the Seidel--Smith construction. 
In addition, we show that the
last theory is equivalent to those of Khovanov. 
More precisely, for a given diagram $D$, let us denote by $\ll D\rr$ 
its Khovanov
hypercube of resolutions. Applying the Khovanov TQFT, we get
a complex $F_{\Kh}\ll D\rr$. On the other hand, 
using the nested Frobenius system,
defined in Section \ref{NFS}, we get the complex $F\ll D\rr$.


\begin{Thm}\label{main} Given a diagram $D$ of a link $L$,
the complexes $F_{\Kh} \ll D\rr$ and $F\ll D\rr$
are isomorphic.
\end{Thm}

Once the equivalence between the geometric  construction of Seidel--Smith 
and the algebraic
one of Cautis--Kamnitzer is established rigorously,
this theorem can be used to finalize the proof of the Seidel--Smith conjecture.

A similar result was independently proved
by a student of C. Stroppel.

The last result of the paper is the
classification of all rank two strict
$\Z_2$--extensions of $\NC$.

\begin{Thm}\label{main1}
 Any strict EQFT based on a $\Z_2$--extension of $\NC$, which agrees with
Khovanov's TQFT after reducing the coefficients modulo 2,
  gives rises to a link invariant homology theory isomorphic
to those of Khovanov.
\end{Thm}

A challenging open problem is to classify all semistrict 
EQFTs based on $\NC$, which associate to a circle a rank two module.
More generally, the problem is to compute the second cohomology of $\NC$
and construct cocycles restricting to the
Schur cocycle of the symmetric group.

An interesting algebraic system underlying
 the categorification of the Kauffman skein module \cite{APS}, \cite{TT}
was proposed recently by Carter and  Saito \cite{CS}. We wonder whether our 
approach could be extended to include their setting.

The paper is organized as follows.
In the first sections we define the categories $\Cob$, $\NC$ and
their extensions. Theorems \ref{main} and \ref{main1} are proved in Section 3.
In the last section, odd Khovanov homology is realized as an extension of $\Cob$.




\subsection*{Acknowledgment} The authors would like to 
thank Catharina Stroppel, Krzysztof Putyra, Alexander Shumakovitch, Christian
Blanchet
and Aaron Lauda
for interesting discussions and to
Dror Bar--Natan for the permission to use his picture of the Khovanov hypercube.

\section{The category of 2--cobordisms and its extensions}

\subsection{The category $\Cob$}
\begin{dfn}{\rm The objects of $\Cob$ are finite ordered
 set of circles.
 The morphisms are isotopy classes of smooth 2--dimensional
cobordisms. The composition is given by gluing of cobordisms.
}\end{dfn}

The category $\Cob$ is a strict symmetric monoidal category with the monoidal
 product given by the ordered disjoint union and the identity
given by the cylinder cobordism.
In particular, we obtain a natural embedding of the symmetric group
in $n$ letters
into the automorphism group of  $n$ circles.

By using Morse theory, one can decompose any 2--cobordism into
 pairs of pants, caps, cups and permutations, proving  the following
well--known presentation of $\Cob$ (see e.g. \cite{H})
\begin{thm}\label{rel-cob}
The morphisms of $\Cob$ are generated by
\begin{center}
\begin{picture}(0,0)%
\epsfig{file=gencob.pstex}%
\end{picture}%
\setlength{\unitlength}{789sp}%
\begingroup\makeatletter\ifx\SetFigFont\undefined%
\gdef\SetFigFont#1#2#3#4#5{%
  \reset@font\fontsize{#1}{#2pt}%
  \fontfamily{#3}\fontseries{#4}\fontshape{#5}%
  \selectfont}%
\fi\endgroup%
\begin{picture}(19012,4860)(7477,-3622)
\put(18001,-3361){\makebox(0,0)[lb]{\smash{{\SetFigFont{12}{14.4}{\familydefault}{\mddefault}{\updefault}$death$}}}}
\put(7651,-3361){\makebox(0,0)[lb]{\smash{{\SetFigFont{12}{14.4}{\familydefault}{\mddefault}{\updefault}$merge$}}}}
\put(11476,-3361){\makebox(0,0)[lb]{\smash{{\SetFigFont{12}{14.4}{\familydefault}{\mddefault}{\updefault}$birth$}}}}
\put(14851,-3361){\makebox(0,0)[lb]{\smash{{\SetFigFont{12}{14.4}{\familydefault}{\mddefault}{\updefault}$split$}}}}
\put(21151,-3361){\makebox(0,0)[lb]{\smash{{\SetFigFont{12}{14.4}{\familydefault}{\mddefault}{\updefault}$permutation$}}}}
\end{picture}%

\end{center}
subject to the following relations:

(1) Commutativity and co-commutativity relation
\begin{center}
\begin{picture}(0,0)%
\epsfig{file=com.pstex}%
\end{picture}%
\setlength{\unitlength}{789sp}%
\begingroup\makeatletter\ifx\SetFigFont\undefined%
\gdef\SetFigFont#1#2#3#4#5{%
  \reset@font\fontsize{#1}{#2pt}%
  \fontfamily{#3}\fontseries{#4}\fontshape{#5}%
  \selectfont}%
\fi\endgroup%
\begin{picture}(11073,3273)(3427,-9010)
\put(9676,-7711){\makebox(0,0)[lb]{\smash{{\SetFigFont{12}{14.4}{\familydefault}{\mddefault}{\updefault}$=$}}}}
\end{picture}%
\hspace{2cm}
\begin{picture}(0,0)%
\epsfig{file=cocom.pstex}%
\end{picture}%
\setlength{\unitlength}{789sp}%
\begingroup\makeatletter\ifx\SetFigFont\undefined%
\gdef\SetFigFont#1#2#3#4#5{%
  \reset@font\fontsize{#1}{#2pt}%
  \fontfamily{#3}\fontseries{#4}\fontshape{#5}%
  \selectfont}%
\fi\endgroup%
\begin{picture}(11073,3273)(14602,-9010)
\put(19426,-7711){\makebox(0,0)[rb]{\smash{{\SetFigFont{12}{14.4}{\familydefault}{\mddefault}{\updefault}$=$}}}}
\end{picture}%

\end{center}

(2) Associativity and coassociativity relations
\begin{center}
\begin{picture}(0,0)%
\epsfig{file=ass.pstex}%
\end{picture}%
\setlength{\unitlength}{789sp}%
\begingroup\makeatletter\ifx\SetFigFont\undefined%
\gdef\SetFigFont#1#2#3#4#5{%
  \reset@font\fontsize{#1}{#2pt}%
  \fontfamily{#3}\fontseries{#4}\fontshape{#5}%
  \selectfont}%
\fi\endgroup%
\begin{picture}(14748,5148)(18277,-10135)
\put(26326,-7861){\makebox(0,0)[rb]{\smash{{\SetFigFont{12}{14.4}{\familydefault}{\mddefault}{\updefault}$=$}}}}
\end{picture}%
\hspace{2cm}
\begin{picture}(0,0)%
\epsfig{file=coass.pstex}%
\end{picture}%
\setlength{\unitlength}{789sp}%
\begingroup\makeatletter\ifx\SetFigFont\undefined%
\gdef\SetFigFont#1#2#3#4#5{%
  \reset@font\fontsize{#1}{#2pt}%
  \fontfamily{#3}\fontseries{#4}\fontshape{#5}%
  \selectfont}%
\fi\endgroup%
\begin{picture}(14748,5148)(3427,-10135)
\put(10126,-7861){\rotatebox{360.0}{\makebox(0,0)[lb]{\smash{{\SetFigFont{12}{14.4}{\familydefault}{\mddefault}{\updefault}$=$}}}}}
\end{picture}%

\end{center}

(3) Frobenius relations
\begin{center}
\begin{picture}(0,0)%
\epsfig{file=frob.pstex}%
\end{picture}%
\setlength{\unitlength}{789sp}%
\begingroup\makeatletter\ifx\SetFigFont\undefined%
\gdef\SetFigFont#1#2#3#4#5{%
  \reset@font\fontsize{#1}{#2pt}%
  \fontfamily{#3}\fontseries{#4}\fontshape{#5}%
  \selectfont}%
\fi\endgroup%
\begin{picture}(21348,5165)(18277,-13068)
\put(25801,-10786){\makebox(0,0)[rb]{\smash{{\SetFigFont{12}{14.4}{\familydefault}{\mddefault}{\updefault}$=$}}}}
\put(33526,-10786){\makebox(0,0)[rb]{\smash{{\SetFigFont{12}{14.4}{\familydefault}{\mddefault}{\updefault}$=$}}}}
\end{picture}%

\end{center}


(4) Unit and Counit relations
\begin{center}
\begin{picture}(0,0)%
\epsfig{file=unitlaw.pstex}%
\end{picture}%
\setlength{\unitlength}{789sp}%
\begingroup\makeatletter\ifx\SetFigFont\undefined%
\gdef\SetFigFont#1#2#3#4#5{%
  \reset@font\fontsize{#1}{#2pt}%
  \fontfamily{#3}\fontseries{#4}\fontshape{#5}%
  \selectfont}%
\fi\endgroup%
\begin{picture}(11448,3222)(3577,-8569)
\put(9901,-7261){\makebox(0,0)[lb]{\smash{{\SetFigFont{12}{14.4}{\familydefault}{\mddefault}{\updefault}$=$}}}}
\end{picture}%
\hspace{2cm}
\begin{picture}(0,0)%
\epsfig{file=counitlaw.pstex}%
\end{picture}%
\setlength{\unitlength}{789sp}%
\begingroup\makeatletter\ifx\SetFigFont\undefined%
\gdef\SetFigFont#1#2#3#4#5{%
  \reset@font\fontsize{#1}{#2pt}%
  \fontfamily{#3}\fontseries{#4}\fontshape{#5}%
  \selectfont}%
\fi\endgroup%
\begin{picture}(11448,3222)(15127,-8569)
\put(20251,-7261){\makebox(0,0)[rb]{\smash{{\SetFigFont{12}{14.4}{\familydefault}{\mddefault}{\updefault}$=$}}}}
\end{picture}%

\end{center}

(5) Permutation relations
\begin{center}
\begin{picture}(0,0)%
\epsfig{file=trans2.pstex}%
\end{picture}%
\setlength{\unitlength}{789sp}%
\begingroup\makeatletter\ifx\SetFigFont\undefined%
\gdef\SetFigFont#1#2#3#4#5{%
  \reset@font\fontsize{#1}{#2pt}%
  \fontfamily{#3}\fontseries{#4}\fontshape{#5}%
  \selectfont}%
\fi\endgroup%
\begin{picture}(10698,3216)(7102,-8944)
\put(13126,-7561){\makebox(0,0)[lb]{\smash{{\SetFigFont{12}{14.4}{\familydefault}{\mddefault}{\updefault}$=$}}}}
\end{picture}%
\hspace{2cm}
\begin{picture}(0,0)%
\epsfig{file=R3trans.pstex}%
\end{picture}%
\setlength{\unitlength}{789sp}%
\begingroup\makeatletter\ifx\SetFigFont\undefined%
\gdef\SetFigFont#1#2#3#4#5{%
  \reset@font\fontsize{#1}{#2pt}%
  \fontfamily{#3}\fontseries{#4}\fontshape{#5}%
  \selectfont}%
\fi\endgroup%
\begin{picture}(19698,5166)(3952,-10294)
\put(12976,-8161){\makebox(0,0)[lb]{\smash{{\SetFigFont{12}{14.4}{\familydefault}{\mddefault}{\updefault}$=$}}}}
\end{picture}%

\end{center}

(6) Unit-Permutations and Counit-Permutation relations
\begin{center}
\begin{picture}(0,0)%
\epsfig{file=permunit.pstex}%
\end{picture}%
\setlength{\unitlength}{789sp}%
\begingroup\makeatletter\ifx\SetFigFont\undefined%
\gdef\SetFigFont#1#2#3#4#5{%
  \reset@font\fontsize{#1}{#2pt}%
  \fontfamily{#3}\fontseries{#4}\fontshape{#5}%
  \selectfont}%
\fi\endgroup%
\begin{picture}(11073,3222)(4027,-7975)
\put(10201,-6661){\makebox(0,0)[lb]{\smash{{\SetFigFont{12}{14.4}{\familydefault}{\mddefault}{\updefault}$=$}}}}
\end{picture}%
\hspace{2cm}
\begin{picture}(0,0)%
\epsfig{file=permcounit.pstex}%
\end{picture}%
\setlength{\unitlength}{789sp}%
\begingroup\makeatletter\ifx\SetFigFont\undefined%
\gdef\SetFigFont#1#2#3#4#5{%
  \reset@font\fontsize{#1}{#2pt}%
  \fontfamily{#3}\fontseries{#4}\fontshape{#5}%
  \selectfont}%
\fi\endgroup%
\begin{picture}(11523,3222)(6952,-8200)
\put(14626,-6736){\makebox(0,0)[rb]{\smash{{\SetFigFont{12}{14.4}{\familydefault}{\mddefault}{\updefault}$=$}}}}
\end{picture}%

\end{center}

(7) Merge-Permutation and Split-Permutation relations
\begin{center}
\begin{picture}(0,0)%
\epsfig{file=permerge.pstex}%
\end{picture}%
\setlength{\unitlength}{789sp}%
\begingroup\makeatletter\ifx\SetFigFont\undefined%
\gdef\SetFigFont#1#2#3#4#5{%
  \reset@font\fontsize{#1}{#2pt}%
  \fontfamily{#3}\fontseries{#4}\fontshape{#5}%
  \selectfont}%
\fi\endgroup%
\begin{picture}(16698,5166)(18277,-10144)
\put(26101,-7861){\makebox(0,0)[rb]{\smash{{\SetFigFont{12}{14.4}{\familydefault}{\mddefault}{\updefault}$=$}}}}
\end{picture}%
\hspace{0.5cm}
\begin{picture}(0,0)%
\epsfig{file=permsplit.pstex}%
\end{picture}%
\setlength{\unitlength}{789sp}%
\begingroup\makeatletter\ifx\SetFigFont\undefined%
\gdef\SetFigFont#1#2#3#4#5{%
  \reset@font\fontsize{#1}{#2pt}%
  \fontfamily{#3}\fontseries{#4}\fontshape{#5}%
  \selectfont}%
\fi\endgroup%
\begin{picture}(16398,5166)(-6548,-10144)
\put(-449,-7861){\rotatebox{360.0}{\makebox(0,0)[lb]{\smash{{\SetFigFont{12}{14.4}{\familydefault}{\mddefault}{\updefault}$=$}}}}}
\end{picture}%

\end{center}

\end{thm}

For  a commutative unital ring $R$, let
 $R$-$\mod$ be the category  of finite projective modules over $R$.
A (1+1)--dimensional topological quantum field theory
(TQFT) is a 
 symmetric (strict) monoidal functor  from $\Cob$ 
to $R$-$\mod$.
 Such TQFTs are in $1:1$ correspondence 
with so--called Frobenius systems
(compare \cite{Ko}). 

One important application of
Frobenius systems is Khovanov's categorification of the  
Jones polynomial \cite{Kho}.

In what follows we will assume that $\Cob$ is a
pre--additive category. This means we  
supply the   set of morphisms (between any two given objects)
with the structure of an abelian group by allowing
 formal $\Z$--linear
combinations of cobordisms and  extend the  composition maps
  bilinearly.

\subsection{Extensions of categories}
In this section we use the language  of 2--categories.
A 2--category is   a category where
 any set of morphisms
has a structure of a category, i.e. we allow morphisms between morphisms
called 2--morphisms.
Given  a 2--category $M$,
the 2--morphisms of $M$ can be composed in two ways.
For any three objects $a,b, c$ of $M$,
the composition in the category $\Mor_{M}(a,b)$ is called 
vertical composition and
the bifunctor $*:\Mor_M(a,b) \times \Mor_M(b,c) \to \Mor_M(a,c)$ is called 
horizontal composition.
These compositions are required to be associative and satisfy an 
interchange law (see
\cite{MacL}  for more details). 


Semistrict monoidal 2--categories can be considered as a weakening of 
monoidal 2--categories, where 
monoidal and
interchange rules hold up to  natural isomorphisms
(compare \cite[Proposition 17]{Lauda}).

Assume $\cal C$ is a strict monoidal category, whose set of morphisms is 
given by
generators and relations.

\begin{dfn}{\rm An extension of $\cal C$ is the semistrict
 monoidal 2--category $\BiC$, which
 has the same set of  
 objects as $\cal C$. The 1--morphisms
of $\BiC$ are compositions of generators of $\cal C$. 
 The 2--morphisms are 
\begin{itemize}
\item
the identity automorphism of any 1--morphism of $\cal C$;
\item
a 2--isomorphism between any
two 1--morphisms subject to a relation in $\C$.
\end{itemize}
}\end{dfn}

This imposes a so--called ``cocycle'' 
condition on the set of 2--morphisms, since
 any composition of 2--morphisms going from a given 1--morphism to itself
should be equal to the identity or any closed loop of 2--morphisms is trivial.

An example of an extension is given by a  weak monoidal category
$(M,\otimes, 1, \alpha, \lambda,\rho)$ where $\alpha$, $\lambda$ and $\rho$
are considered as 2--isomorphisms and 
the cocycle condition holds due
  to MacLane's coherence theorem \cite[Chapter VII]{MacL}.
In the case, when $\cal C$ is $\Cob$ restricted to connected cobordisms,
(i.e. permutation is removed from the set of generators in 
Theorem \ref{rel-cob}
as well as relations (1),(5), (6) and (7)), then any 
 pseudo Frobenius algebra, described in 
\cite[Proposition 25]{Lauda}, defines an extension $\BiC$. The cocycle
condition holds due to Lemmas 32, 33 in \cite{Lauda}.

Providing $\cal C$ with a structure of a pre--additive category, 
we have a natural 
$\Z$--action on the set of 1--morphisms, restricting to $\Z_2=\{1,-1\}$,
the group of two elements written multiplicatively, we can define
a $\Z_2$--extension of $\cal C$, in which the 2--morphisms are just 
plus or minus the identity. Note that in this case $\BiC$
can be considered as a weak monoidal category, with the same set of 
generating
1--morphisms as $\cal C$, but with sign modified relations.

For  any cobordism category $\cal C$, an {\it extended
quantum field theory} (EQFT) based on $\cal C$ is 
a bifunctor from $\BiC$
 to $R$-$\mod$, mapping 2--morphisms  to
natural transformations of $R$--modules. The EQFT is called strict
if $\BiC$ is strict.

\section{Embedded cobordisms }

Let $S^{\amalg a}$ be the disjoint union of $a$ copies of a circle smoothly 
embedded into a plane. Note that the embedding induces a partial
order on the set of circles as follows. For two  circles $c_1$ and $c_2$,
we say $c_1< c_2$, 
if $c_1$ is inside $c_2$. 

\begin{dfn}{\rm The objects of $\NC$ are finite collections of circles
embedded into a plane. The morphisms are generated by
\begin{center}
\begin{picture}(0,0)%
\epsfig{file=gencob.pstex}%
\end{picture}%
\setlength{\unitlength}{789sp}%
\begingroup\makeatletter\ifx\SetFigFont\undefined%
\gdef\SetFigFont#1#2#3#4#5{%
  \reset@font\fontsize{#1}{#2pt}%
  \fontfamily{#3}\fontseries{#4}\fontshape{#5}%
  \selectfont}%
\fi\endgroup%
\begin{picture}(19012,4860)(7477,-3622)
\put(18001,-3361){\makebox(0,0)[lb]{\smash{{\SetFigFont{12}{14.4}{\familydefault}{\mddefault}{\updefault}$death$}}}}
\put(7651,-3361){\makebox(0,0)[lb]{\smash{{\SetFigFont{12}{14.4}{\familydefault}{\mddefault}{\updefault}$merge$}}}}
\put(11476,-3361){\makebox(0,0)[lb]{\smash{{\SetFigFont{12}{14.4}{\familydefault}{\mddefault}{\updefault}$birth$}}}}
\put(14851,-3361){\makebox(0,0)[lb]{\smash{{\SetFigFont{12}{14.4}{\familydefault}{\mddefault}{\updefault}$split$}}}}
\put(21151,-3361){\makebox(0,0)[lb]{\smash{{\SetFigFont{12}{14.4}{\familydefault}{\mddefault}{\updefault}$permutation$}}}}
\end{picture}%
\\
\vspace{1cm}
\begin{picture}(0,0)%
\epsfig{file=genestedcob.pstex}%
\end{picture}%
\setlength{\unitlength}{789sp}%
\begingroup\makeatletter\ifx\SetFigFont\undefined%
\gdef\SetFigFont#1#2#3#4#5{%
  \reset@font\fontsize{#1}{#2pt}%
  \fontfamily{#3}\fontseries{#4}\fontshape{#5}%
  \selectfont}%
\fi\endgroup%
\begin{picture}(10373,4190)(24076,-9022)
\put(30001,-8761){\makebox(0,0)[lb]{\smash{{\SetFigFont{12}{14.4}{\familydefault}{\mddefault}{\updefault}$nested\mbox{ }merge$}}}}
\put(24076,-8761){\makebox(0,0)[lb]{\smash{{\SetFigFont{12}{14.4}{\familydefault}{\mddefault}{\updefault}$nested\mbox{ }split$}}}}
\end{picture}%

\end{center}

subject to the 
following sets of relations:

(1) Frobenius type relations

\begin{center}
\begin{picture}(0,0)%
\epsfig{file=1.pstex}%
\end{picture}%
\setlength{\unitlength}{789sp}%
\begingroup\makeatletter\ifx\SetFigFont\undefined%
\gdef\SetFigFont#1#2#3#4#5{%
  \reset@font\fontsize{#1}{#2pt}%
  \fontfamily{#3}\fontseries{#4}\fontshape{#5}%
  \selectfont}%
\fi\endgroup%
\begin{picture}(17728,4546)(-41127,-8484)
\put(-33824,-6511){\makebox(0,0)[lb]{\smash{{\SetFigFont{12}{14.4}{\familydefault}{\mddefault}{\updefault}$=$}}}}
\end{picture}%
\\
\begin{picture}(0,0)%
\epsfig{file=2.pstex}%
\end{picture}%
\setlength{\unitlength}{789sp}%
\begingroup\makeatletter\ifx\SetFigFont\undefined%
\gdef\SetFigFont#1#2#3#4#5{%
  \reset@font\fontsize{#1}{#2pt}%
  \fontfamily{#3}\fontseries{#4}\fontshape{#5}%
  \selectfont}%
\fi\endgroup%
\begin{picture}(17728,4546)(-23398,-8484)
\put(-12974,-6511){\rotatebox{360.0}{\makebox(0,0)[rb]{\smash{{\SetFigFont{12}{14.4}{\familydefault}{\mddefault}{\updefault}$=$}}}}}
\end{picture}%

\begin{picture}(0,0)%
\epsfig{file=11.pstex}%
\end{picture}%
\setlength{\unitlength}{789sp}%
\begingroup\makeatletter\ifx\SetFigFont\undefined%
\gdef\SetFigFont#1#2#3#4#5{%
  \reset@font\fontsize{#1}{#2pt}%
  \fontfamily{#3}\fontseries{#4}\fontshape{#5}%
  \selectfont}%
\fi\endgroup%
\begin{picture}(17728,4546)(-41127,-8484)
\put(-33824,-6511){\rotatebox{360.0}{\makebox(0,0)[lb]{\smash{{\SetFigFont{12}{14.4}{\familydefault}{\mddefault}{\updefault}$=$}}}}}
\end{picture}%
\\
 \vspace{0.3cm}
\begin{picture}(0,0)%
\epsfig{file=5.pstex}%
\end{picture}%
\setlength{\unitlength}{789sp}%
\begingroup\makeatletter\ifx\SetFigFont\undefined%
\gdef\SetFigFont#1#2#3#4#5{%
  \reset@font\fontsize{#1}{#2pt}%
  \fontfamily{#3}\fontseries{#4}\fontshape{#5}%
  \selectfont}%
\fi\endgroup%
\begin{picture}(14783,4935)(9242,-14484)
\put(16201,-12286){\makebox(0,0)[lb]{\smash{{\SetFigFont{12}{14.4}{\familydefault}{\mddefault}{\updefault}$=$}}}}
\end{picture}%
\hspace{1cm}
\begin{picture}(0,0)%
\epsfig{file=9.pstex}%
\end{picture}%
\setlength{\unitlength}{789sp}%
\begingroup\makeatletter\ifx\SetFigFont\undefined%
\gdef\SetFigFont#1#2#3#4#5{%
  \reset@font\fontsize{#1}{#2pt}%
  \fontfamily{#3}\fontseries{#4}\fontshape{#5}%
  \selectfont}%
\fi\endgroup%
\begin{picture}(14782,4935)(-5873,-14484)
\put(1951,-12286){\rotatebox{360.0}{\makebox(0,0)[rb]{\smash{{\SetFigFont{12}{14.4}{\familydefault}{\mddefault}{\updefault}$=$}}}}}
\end{picture}%
\\
 \vspace{0.3cm}
\begin{picture}(0,0)%
\epsfig{file=3.pstex}%
\end{picture}%
\setlength{\unitlength}{789sp}%
\begingroup\makeatletter\ifx\SetFigFont\undefined%
\gdef\SetFigFont#1#2#3#4#5{%
  \reset@font\fontsize{#1}{#2pt}%
  \fontfamily{#3}\fontseries{#4}\fontshape{#5}%
  \selectfont}%
\fi\endgroup%
\begin{picture}(16073,4402)(-22198,-8484)
\put(-12974,-6661){\rotatebox{360.0}{\makebox(0,0)[rb]{\smash{{\SetFigFont{12}{14.4}{\familydefault}{\mddefault}{\updefault}$=$}}}}}
\end{picture}%
\hspace{1cm}
\begin{picture}(0,0)%
\epsfig{file=8.pstex}%
\end{picture}%
\setlength{\unitlength}{789sp}%
\begingroup\makeatletter\ifx\SetFigFont\undefined%
\gdef\SetFigFont#1#2#3#4#5{%
  \reset@font\fontsize{#1}{#2pt}%
  \fontfamily{#3}\fontseries{#4}\fontshape{#5}%
  \selectfont}%
\fi\endgroup%
\begin{picture}(16073,4402)(-6023,-8484)
\put(826,-6661){\makebox(0,0)[lb]{\smash{{\SetFigFont{12}{14.4}{\familydefault}{\mddefault}{\updefault}$=$}}}}
\end{picture}%

\end{center}

(2) Associativity type relations

\begin{center}
\begin{picture}(0,0)%
\epsfig{file=A.pstex}%
\end{picture}%
\setlength{\unitlength}{789sp}%
\begingroup\makeatletter\ifx\SetFigFont\undefined%
\gdef\SetFigFont#1#2#3#4#5{%
  \reset@font\fontsize{#1}{#2pt}%
  \fontfamily{#3}\fontseries{#4}\fontshape{#5}%
  \selectfont}%
\fi\endgroup%
\begin{picture}(14858,4935)(20491,-9159)
\put(27076,-7186){\makebox(0,0)[lb]{\smash{{\SetFigFont{12}{14.4}{\familydefault}{\mddefault}{\updefault}$=$}}}}
\end{picture}%
\hspace{1cm}
\begin{picture}(0,0)%
\epsfig{file=B.pstex}%
\end{picture}%
\setlength{\unitlength}{789sp}%
\begingroup\makeatletter\ifx\SetFigFont\undefined%
\gdef\SetFigFont#1#2#3#4#5{%
  \reset@font\fontsize{#1}{#2pt}%
  \fontfamily{#3}\fontseries{#4}\fontshape{#5}%
  \selectfont}%
\fi\endgroup%
\begin{picture}(17920,4621)(19454,-13659)
\put(27676,-11761){\makebox(0,0)[lb]{\smash{{\SetFigFont{12}{14.4}{\familydefault}{\mddefault}{\updefault}$=$}}}}
\end{picture}%
\\
\vspace{0.3cm}
\begin{picture}(0,0)%
\epsfig{file=C.pstex}%
\end{picture}%
\setlength{\unitlength}{789sp}%
\begingroup\makeatletter\ifx\SetFigFont\undefined%
\gdef\SetFigFont#1#2#3#4#5{%
  \reset@font\fontsize{#1}{#2pt}%
  \fontfamily{#3}\fontseries{#4}\fontshape{#5}%
  \selectfont}%
\fi\endgroup%
\begin{picture}(17845,4546)(19529,-13584)
\put(27676,-11611){\makebox(0,0)[lb]{\smash{{\SetFigFont{12}{14.4}{\familydefault}{\mddefault}{\updefault}$=$}}}}
\end{picture}%
\hspace{1cm}
\begin{picture}(0,0)%
\epsfig{file=D.pstex}%
\end{picture}%
\setlength{\unitlength}{789sp}%
\begingroup\makeatletter\ifx\SetFigFont\undefined%
\gdef\SetFigFont#1#2#3#4#5{%
  \reset@font\fontsize{#1}{#2pt}%
  \fontfamily{#3}\fontseries{#4}\fontshape{#5}%
  \selectfont}%
\fi\endgroup%
\begin{picture}(15353,4402)(22547,-9534)
\put(30226,-7861){\makebox(0,0)[lb]{\smash{{\SetFigFont{12}{14.4}{\familydefault}{\mddefault}{\updefault}$=$}}}}
\end{picture}%

\end{center}

(3) Coassociativity type relations

\begin{center}
\begin{picture}(0,0)%
\epsfig{file=A1.pstex}%
\end{picture}%
\setlength{\unitlength}{789sp}%
\begingroup\makeatletter\ifx\SetFigFont\undefined%
\gdef\SetFigFont#1#2#3#4#5{%
  \reset@font\fontsize{#1}{#2pt}%
  \fontfamily{#3}\fontseries{#4}\fontshape{#5}%
  \selectfont}%
\fi\endgroup%
\begin{picture}(14858,4935)(5602,-9159)
\put(13876,-7186){\makebox(0,0)[rb]{\smash{{\SetFigFont{12}{14.4}{\familydefault}{\mddefault}{\updefault}$=$}}}}
\end{picture}%
\hspace{1cm}
\begin{picture}(0,0)%
\epsfig{file=B1.pstex}%
\end{picture}%
\setlength{\unitlength}{789sp}%
\begingroup\makeatletter\ifx\SetFigFont\undefined%
\gdef\SetFigFont#1#2#3#4#5{%
  \reset@font\fontsize{#1}{#2pt}%
  \fontfamily{#3}\fontseries{#4}\fontshape{#5}%
  \selectfont}%
\fi\endgroup%
\begin{picture}(17920,4621)(1477,-13659)
\put(11176,-11761){\makebox(0,0)[rb]{\smash{{\SetFigFont{12}{14.4}{\familydefault}{\mddefault}{\updefault}$=$}}}}
\end{picture}%
\\
\vspace{0.3cm}
\begin{picture}(0,0)%
\epsfig{file=C1.pstex}%
\end{picture}%
\setlength{\unitlength}{789sp}%
\begingroup\makeatletter\ifx\SetFigFont\undefined%
\gdef\SetFigFont#1#2#3#4#5{%
  \reset@font\fontsize{#1}{#2pt}%
  \fontfamily{#3}\fontseries{#4}\fontshape{#5}%
  \selectfont}%
\fi\endgroup%
\begin{picture}(17845,4546)(1927,-13584)
\put(11626,-11611){\makebox(0,0)[rb]{\smash{{\SetFigFont{12}{14.4}{\familydefault}{\mddefault}{\updefault}$=$}}}}
\end{picture}%
\hspace{1cm}
\begin{picture}(0,0)%
\epsfig{file=D1.pstex}%
\end{picture}%
\setlength{\unitlength}{789sp}%
\begingroup\makeatletter\ifx\SetFigFont\undefined%
\gdef\SetFigFont#1#2#3#4#5{%
  \reset@font\fontsize{#1}{#2pt}%
  \fontfamily{#3}\fontseries{#4}\fontshape{#5}%
  \selectfont}%
\fi\endgroup%
\begin{picture}(15352,4402)(7102,-9534)
\put(14776,-7861){\makebox(0,0)[rb]{\smash{{\SetFigFont{12}{14.4}{\familydefault}{\mddefault}{\updefault}$=$}}}}
\end{picture}%

\end{center}

(4) Cancellation
\begin{center}
\begin{picture}(0,0)%
\epsfig{file=nestedcounitlaw.pstex}%
\end{picture}%
\setlength{\unitlength}{789sp}%
\begingroup\makeatletter\ifx\SetFigFont\undefined%
\gdef\SetFigFont#1#2#3#4#5{%
  \reset@font\fontsize{#1}{#2pt}%
  \fontfamily{#3}\fontseries{#4}\fontshape{#5}%
  \selectfont}%
\fi\endgroup%
\begin{picture}(11748,2466)(7552,-12244)
\put(14401,-11386){\makebox(0,0)[lb]{\smash{{\SetFigFont{12}{14.4}{\familydefault}{\mddefault}{\updefault}$=$}}}}
\end{picture}%
\hspace{2cm}
\begin{picture}(0,0)%
\epsfig{file=nstedunitlaw.pstex}%
\end{picture}%
\setlength{\unitlength}{789sp}%
\begingroup\makeatletter\ifx\SetFigFont\undefined%
\gdef\SetFigFont#1#2#3#4#5{%
  \reset@font\fontsize{#1}{#2pt}%
  \fontfamily{#3}\fontseries{#4}\fontshape{#5}%
  \selectfont}%
\fi\endgroup%
\begin{picture}(11748,2466)(-4298,-12244)
\put(601,-11386){\rotatebox{360.0}{\makebox(0,0)[rb]{\smash{{\SetFigFont{12}{14.4}{\familydefault}{\mddefault}{\updefault}$=$}}}}}
\end{picture}%

\end{center}

(5) Torus relation 
\begin{center}
\begin{picture}(0,0)%
\epsfig{file=torus.pstex}%
\end{picture}%
\setlength{\unitlength}{789sp}%
\begingroup\makeatletter\ifx\SetFigFont\undefined%
\gdef\SetFigFont#1#2#3#4#5{%
  \reset@font\fontsize{#1}{#2pt}%
  \fontfamily{#3}\fontseries{#4}\fontshape{#5}%
  \selectfont}%
\fi\endgroup%
\begin{picture}(14673,3198)(6502,-14710)
\put(13276,-13486){\makebox(0,0)[lb]{\smash{{\SetFigFont{12}{14.4}{\familydefault}{\mddefault}{\updefault}$=$}}}}
\end{picture}%

\end{center}

}
In addition, the merge, the split, the birth, the death and the permutation are still subject to all the relations of Theorem \ref{rel-cob}.
\label{nc}
\end{dfn}

\begin{remark}
The relations of Definition \ref{nc} can also be described as follows:
\begin{center}
\begin{picture}(0,0)%
\epsfig{file=assface.pstex}%
\end{picture}%
\setlength{\unitlength}{789sp}%
\begingroup\makeatletter\ifx\SetFigFont\undefined%
\gdef\SetFigFont#1#2#3#4#5{%
  \reset@font\fontsize{#1}{#2pt}%
  \fontfamily{#3}\fontseries{#4}\fontshape{#5}%
  \selectfont}%
\fi\endgroup%
\begin{picture}(19171,3914)(7478,-7976)
\end{picture}%
\\
\begin{picture}(0,0)%
\epsfig{file=coassface.pstex}%
\end{picture}%
\setlength{\unitlength}{789sp}%
\begingroup\makeatletter\ifx\SetFigFont\undefined%
\gdef\SetFigFont#1#2#3#4#5{%
  \reset@font\fontsize{#1}{#2pt}%
  \fontfamily{#3}\fontseries{#4}\fontshape{#5}%
  \selectfont}%
\fi\endgroup%
\begin{picture}(32622,5072)(2802,-9197)
\end{picture}%
\\
\begin{picture}(0,0)%
\epsfig{file=frobface.pstex}%
\end{picture}%
\setlength{\unitlength}{789sp}%
\begingroup\makeatletter\ifx\SetFigFont\undefined%
\gdef\SetFigFont#1#2#3#4#5{%
  \reset@font\fontsize{#1}{#2pt}%
  \fontfamily{#3}\fontseries{#4}\fontshape{#5}%
  \selectfont}%
\fi\endgroup%
\begin{picture}(19687,7853)(3245,-12160)
\end{picture}%
\\
\begin{picture}(0,0)%
\epsfig{file=torusface.pstex}%
\end{picture}%
\setlength{\unitlength}{789sp}%
\begingroup\makeatletter\ifx\SetFigFont\undefined%
\gdef\SetFigFont#1#2#3#4#5{%
  \reset@font\fontsize{#1}{#2pt}%
  \fontfamily{#3}\fontseries{#4}\fontshape{#5}%
  \selectfont}%
\fi\endgroup%
\begin{picture}(2012,2458)(3838,-11269)
\end{picture}%

\end{center}
where the black circle corresponds to the starting configuration of circles, and the dashed arcs correspond to the operations
 which are performed. Notice that  changing the order of 
operations produce the two different sides of the relations in Definition \ref{nc}. In addition, associativity of the merge, the coassociativity of the split and the usual Frobenius relation are also depicted here.
\end{remark}

The category $\NC$ is a symmetric strict monoidal category with a 
tensor product given by a partially ordered disjoint union,
i.e. circles on the same level of nestedness are ordered.
In particular, we obtain a natural embedding of the symmetric group
into the automorphisms of any object, permuting 
circles not ordered by nestedness and at the same level of nestedness.\\

Any morphism in $\NC$ is the composite of such a permutation
and the tensor product of connected morphisms of $\NC$.


\begin{lem}
Any connected morphism in $\NC$ has the following normal form:
\begin{center}
\begin{picture}(0,0)%
\epsfig{file=normalformembedded.pstex}%
\end{picture}%
\setlength{\unitlength}{671sp}%
\begingroup\makeatletter\ifx\SetFigFont\undefined%
\gdef\SetFigFont#1#2#3#4#5{%
  \reset@font\fontsize{#1}{#2pt}%
  \fontfamily{#3}\fontseries{#4}\fontshape{#5}%
  \selectfont}%
\fi\endgroup%
\begin{picture}(24325,9369)(6704,-18032)
\end{picture}%

\end{center}
\end{lem}

\begin{proof}
Assume that the boundary of our connected genus $g$ cobordism  $C$
consists of $a$ incoming circles and $b$ outgoing ones.
Let us suppose that $C$ is a composition of
 $B$ births, $D$ deaths, $M$ merges and $S$ splits.
Then we have
$$2-2g-a-b=B+D-M-S, \quad \quad a-M+B=b-S+D$$
or
$$ M=a+g-1+B,\quad \quad S=b+g-1+D$$
We arrive at the normal form if we will be able to push all
merges (resp. splits) to the incoming (resp. outgoing) boundary of $C$.
 From the above formulas
we see that
$B$ merges and $D$ splits will cancel with the births and deaths, respectively,
and $g$ splits and merges put together will create $g$ handles.
The remaining $a-1$  merges 
commute  with any split (nested or not nested one)
due to the  Frobenius type relations. Finally using the associativity type
relations, we can commute nested and unnested merges (resp. splits) 
 and 
 arrive at the form in Figure \ref{inhandle}.

\begin{figure}[!h]
\begin{center}
\input{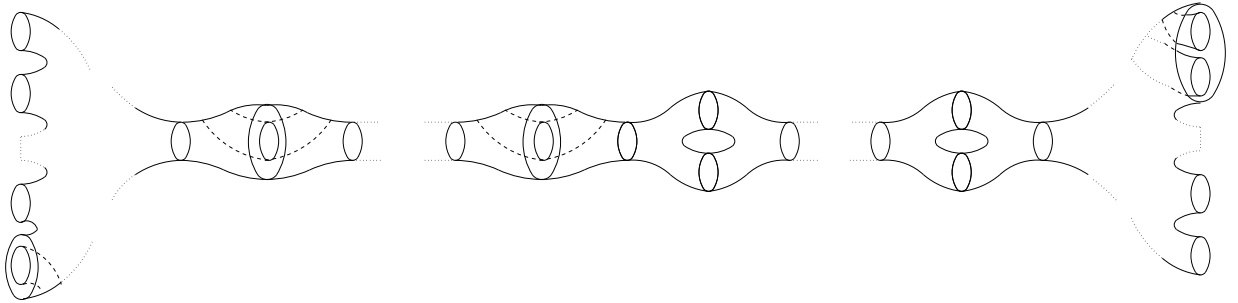}
\end{center}
\caption{Normal form with inner and outer handles}
\label{inhandle}
\end{figure}

 Furthermore applying the Torus relation, one can now reduce to the 
normal form.
\end{proof}

\subsection{Embedded cobordisms}

\begin{dfn}{\rm
A smoothly embedded  2--dimensional cobordism  from $S^{\amalg a}$
to $S^{\amalg b}$ is a pair $(F,\phi)$, where
$F$ a smooth 2--dimensional surface  whose
boundary consists of $a+b$ circles  and
$\phi: F \hookrightarrow \R^2 \times [0,1]$ is a smooth embedding, such that
$\phi|_{\partial F} \cap \R^2\times \{0\}= S^{\amalg a}$ and
$\phi|_{\partial F} \cap \R^2\times \{1\}= S^{\amalg b}$.
 }

\end{dfn}
\begin{dfn}\label{ec} {\rm
The objects of $\EC$ are circles smoothly embedded into a plane.
The morphisms are isotopy   classes of smoothly embedded 2--dimensional 
cobordisms subject to the unknotting relation:
\begin{center}
\begin{picture}(0,0)%
\epsfig{file=switchembedded.pstex}%
\end{picture}%
\setlength{\unitlength}{789sp}%
\begingroup\makeatletter\ifx\SetFigFont\undefined%
\gdef\SetFigFont#1#2#3#4#5{%
  \reset@font\fontsize{#1}{#2pt}%
  \fontfamily{#3}\fontseries{#4}\fontshape{#5}%
  \selectfont}%
\fi\endgroup%
\begin{picture}(19016,8104)(16668,-6963)
\put(25351,-3136){\makebox(0,0)[lb]{\smash{{\SetFigFont{12}{14.4}{\familydefault}{\mddefault}{\updefault}$=$}}}}
\end{picture}%

\end{center}
The composition is given by gluing along the boundary.

}
\end{dfn}

The category $\EC$ is again a symmetric strict monoidal category with
a tensor product given by the partially ordered disjoint union and with the
action of the permutation group depending on nestedness.

\begin{thm} The category $\EC$ is isomorphic to
the category $\NC$.
\end{thm}

\begin{proof}
By \cite{H}, any smooth 2--cobordism allows a pair of pants
decomposition. Modulo the unknotting relation,
there are two ways to embed a pair
of pants into $\R^3$, providing 
the list of generators
 in Definition \ref{nc}. 
The relations do not change the isotopy class
of an embedded cobordism and allow to bring it into a normal form.
It remains just to say, that the normal forms of two equivalent connected
 cobordisms
coincide. 
\end{proof}


\subsection{Strict $\Z_2$--extension}
 Assume $\NC$ is a
pre--additive category. 



\begin{dfn}
{\rm Let $\NC_1$ be the strict monoidal 2--category obtained 
from $\NC$ by replacing
the torus relation with
\begin{center}
(T1)\begin{picture}(0,0)%
\epsfig{file=torusminus.pstex}%
\end{picture}%
\setlength{\unitlength}{789sp}%
\begingroup\makeatletter\ifx\SetFigFont\undefined%
\gdef\SetFigFont#1#2#3#4#5{%
  \reset@font\fontsize{#1}{#2pt}%
  \fontfamily{#3}\fontseries{#4}\fontshape{#5}%
  \selectfont}%
\fi\endgroup%
\begin{picture}(14673,3198)(6502,-14710)
\put(13276,-13486){\makebox(0,0)[lb]{\smash{{\SetFigFont{12}{14.4}{\familydefault}{\mddefault}{\updefault}$=-$}}}}
\end{picture}%

\end{center}}

\end{dfn}

\begin{lem}
$\NC_1$ is a $\Z_2$--extension of $\NC$.
\label{cocy}
\end{lem}

\begin{proof}
The only non--trivial 2--morphism corresponds to
the torus relation. It remains to show that the automorphism
group of any 1--morphism is trivial. By 
Bergman's Diamond Lemma \cite{Berg}, 
 it suffices to check
that any cube with  T1 face has an even number of anticommutative faces. 
This is a simple case by case check. The 10 cubes to check are 
depicted in Figure \ref{cocycle}.
\begin{figure}[!h]
\begin{center}
\input{cocyclenested.pstex_t}
\end{center}
\caption{Cocycle conditions for $\NC_1$}
\label{cocycle}
\end{figure}

\end{proof}

\begin{remark}
Notice that any element of $\NC_1$ does still have a normal form, which corresponds to the usual one
 plus the  information of the parity of the  number of inner 1--handles in Figure \ref{inhandle}.
\end{remark}

\subsection{Nested Frobenius system}\label{NFS}
In this section we construct a strict EQFT based on $\NC_1$,
as proposed by Stroppel \cite{C-talk}.

As in \cite{Kho}, let us consider the 2--dimensional module
$A:=\Z[t]\la\one,X\ra$ over the polynomial ring $R:=\Z[t]$ in one
variable. We denote by $\one$ the 
 image of 1 under the embedding $\eta: R\to A$.
For $\e\in \{0,1\}$,
we define two kinds of a multiplication
$m_\e:A\otimes A\to A$ as follows:
\be
m_\e:\left\{\begin{array}{l@{\quad\mapsto\quad}r}
\one\otimes \one & \one\\
\one\otimes X& X\\
X\otimes \one& (-1)^\e X\\
X\otimes X& (-1)^\e t
\end{array}
\right.
\ee
Further, we define two comultiplications
$\Delta_\e: A\to A\otimes A $ and a counit $\epsilon: A\to R$
as follows.
\be
\Delta_\e:\left\{\begin{array}{l@{\quad\mapsto\quad}r}
 \one & X\otimes\one+ (-1)^\e\;  \one\otimes X\\
 X& 
X\otimes X + (-1)^\e \; t \one\otimes \one\\
\end{array}
\right.
\quad\quad
\epsilon:\left\{\begin{array}{l@{\quad\mapsto\quad}r}
 \one & 0\\
 X& 1\\
\end{array}
\right.
\ee
The functor $F: \NC_1\to R$-$\mod$ maps any object
$S^{\amalg a}$ to
 $A^{\otimes a}$ and is defined on the generating morphisms
as follows:
\be 
F\left(\begin{array}{c}\begin{picture}(0,0)%
\epsfig{file=merge.pstex}%
\end{picture}%
\setlength{\unitlength}{395sp}%
\begingroup\makeatletter\ifx\SetFigFont\undefined%
\gdef\SetFigFont#1#2#3#4#5{%
  \reset@font\fontsize{#1}{#2pt}%
  \fontfamily{#3}\fontseries{#4}\fontshape{#5}%
  \selectfont}%
\fi\endgroup%
\begin{picture}(3198,3198)(7477,-2110)
\end{picture}%
\end{array}\right)=m_0, \quad F\left(\begin{array}{c}\begin{picture}(0,0)%
\epsfig{file=split.pstex}%
\end{picture}%
\setlength{\unitlength}{395sp}%
\begingroup\makeatletter\ifx\SetFigFont\undefined%
\gdef\SetFigFont#1#2#3#4#5{%
  \reset@font\fontsize{#1}{#2pt}%
  \fontfamily{#3}\fontseries{#4}\fontshape{#5}%
  \selectfont}%
\fi\endgroup%
\begin{picture}(3198,3198)(14227,-2035)
\end{picture}%
\end{array}\right)=\Delta_0,
\ee
\be 
F\left(\begin{array}{c}\begin{picture}(0,0)%
\epsfig{file=nestedmerge.pstex}%
\end{picture}%
\setlength{\unitlength}{395sp}%
\begingroup\makeatletter\ifx\SetFigFont\undefined%
\gdef\SetFigFont#1#2#3#4#5{%
  \reset@font\fontsize{#1}{#2pt}%
  \fontfamily{#3}\fontseries{#4}\fontshape{#5}%
  \selectfont}%
\fi\endgroup%
\begin{picture}(3652,2412)(30797,-7244)
\end{picture}%
\end{array}\right)=m_1, \quad F\left(\begin{array}{c}\begin{picture}(0,0)%
\epsfig{file=nestedsplit.pstex}%
\end{picture}%
\setlength{\unitlength}{395sp}%
\begingroup\makeatletter\ifx\SetFigFont\undefined%
\gdef\SetFigFont#1#2#3#4#5{%
  \reset@font\fontsize{#1}{#2pt}%
  \fontfamily{#3}\fontseries{#4}\fontshape{#5}%
  \selectfont}%
\fi\endgroup%
\begin{picture}(3652,2412)(24427,-7289)
\end{picture}%
\end{array}\right)=\Delta_1,
\ee
\be 
F\left(\begin{array}{c}\begin{picture}(0,0)%
\epsfig{file=birth.pstex}%
\end{picture}%
\setlength{\unitlength}{395sp}%
\begingroup\makeatletter\ifx\SetFigFont\undefined%
\gdef\SetFigFont#1#2#3#4#5{%
  \reset@font\fontsize{#1}{#2pt}%
  \fontfamily{#3}\fontseries{#4}\fontshape{#5}%
  \selectfont}%
\fi\endgroup%
\begin{picture}(1098,1278)(11977,-1075)
\end{picture}%
\end{array}\right)=\eta, \quad F\left(\begin{array}{c}\begin{picture}(0,0)%
\epsfig{file=death.pstex}%
\end{picture}%
\setlength{\unitlength}{395sp}%
\begingroup\makeatletter\ifx\SetFigFont\undefined%
\gdef\SetFigFont#1#2#3#4#5{%
  \reset@font\fontsize{#1}{#2pt}%
  \fontfamily{#3}\fontseries{#4}\fontshape{#5}%
  \selectfont}%
\fi\endgroup%
\begin{picture}(1098,1278)(18802,-1075)
\end{picture}%
\end{array}\right)=\epsilon,
\ee

 The convention is that 
in $A\otimes A$ the first factor corresponds to the inner
circle. It is easy to see that $F$ preserves all the relations listed in Definition
\ref{nc}.

Let us introduce a grading on $A$ by putting
$$\gr(t):=-4,\quad \gr(X):=-1,\quad \gr(\one):=1$$
On the tensor product $A^{\otimes n}$ the grading is given
by $\gr(a_1\otimes\dots \otimes a_n):=\gr(a_1)+\dots+\gr(a_n)$.

There exist a natural grading on $\NC$ given by the Euler characteristic of
cobordisms. As in Khovanov's case if $t=0$,
$F$ is grading preserving. 

\section{Nested Khovanov homology}

\subsection{Khovanov's hypercube}\label{11}

\begin{figure}
\mbox{\epsfysize=7cm \epsffile{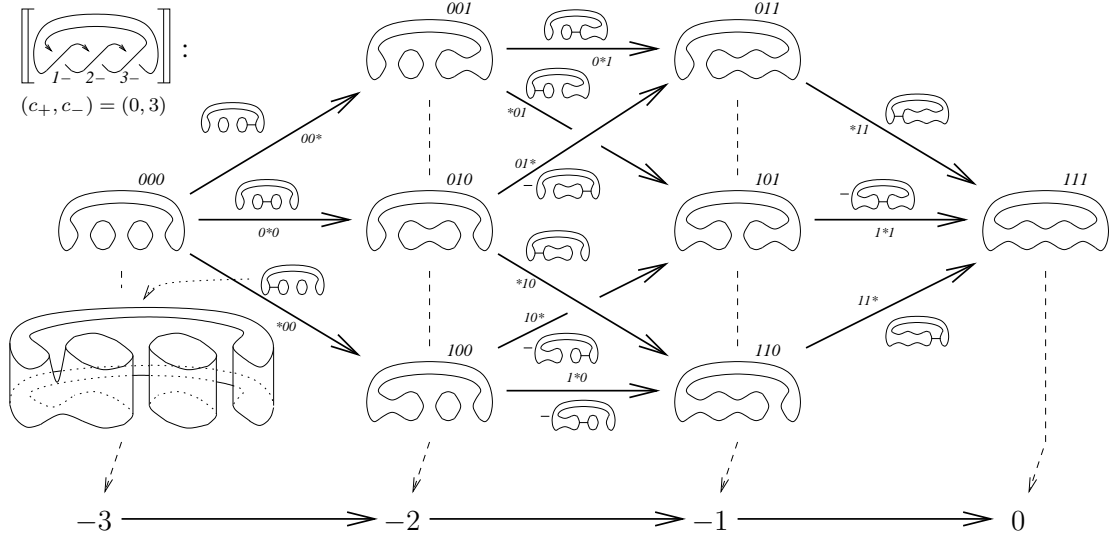}}
\v8
\caption{  The cube of resolutions  for the trefoil}
\label{Diag}
\end{figure}

Suppose we have  a generic diagram $D$ of an oriented link $L$ in $S^3$
with $c$ crossings. By resolving  crossings of $D$ in two ways  as  
prescribed by the Kauffman skein relations, one can
 associate to $D$  a $c$--dimensional cube 
of resolutions (compare \cite{Kho}
or \cite{BN1}).
The vertices of the cube
correspond to the configurations of circles obtained after smoothing
of all crossings in $D$.  For any crossing, two different smoothings are
allowed: the  $0$--  and  the $1$--smoothings. Therefore,
we have $2^c$ vertices. After numbering the crossings of $D$,
we can label the vertices of the cube by $c$--letter
strings of $0$'s and $1$'s,
specifying  the smoothing  chosen at each crossing.
The cube is skewered along its main diagonal,
from $00...0$ to $11...1$.
The number of 1 in the labeling of a vertex is
equal to its `height' $k$. The cube is displayed in
such a way that the vertices of height $k$ project down
to the point $r:=k-c_-$, where $c_\pm$ are the numbers of 
positive, resp. negative crossings in $D$ (see Figure \ref{Diag}).

\v8

Two vertices of the hypercube are connected by an edge
if their labellings  differ by one letter. In Figure \ref{Diag},
 this letter is labeled by $*$.
 The edges are directed
(from the vertex where this letter is $0$ to the vertex where
it is $1$). The edges correspond to a saddle cobordisms from the tail
configuration of circles to the head configuration (compare Figure \ref{Diag}).

We denote
this hypercube of resolutions by $\ll D\rr$,
and would like to interpret it
 as a complex.
The $r$th chain ``space'' $\ll D\rr^{r}$ 
is a  formal direct sum of the
$\frac{c!}{k!(c-k)!}$ ``spaces'' at height $k$
in the hypercube and  the sum of ``maps''
with tails at height $k$
defines the  $r$th differential. 
To achieve $\partial^2=0$, we
assign a minus to any edge which has an odd number of 1's before $*$
in its labeling. 

Applying (1+1) TQFT
 $F_{Kh}$ to $\ll D\rr$, 
which  sends any merge to $m_0$ and any split to $\Delta_0$, we get
a complex of $R$-modules
$(F_{\Kh}\ll D\rr, \partial_{\Kh})$. Its graded homology groups,
known as Khovanov homology, are link invariants and
the graded Euler characteristic is given by the Jones polynomial.

\subsection{Nested homology}
Applying $F$ to the Khovanov hypercube, one can
define a chain complex as follows.
The $r$th chain group will be
$F\ll D\rr^{r}$, the image of $\ll D\rr^{r}$ after applying 
the functor $F$, and the maps are defined by applying $F$ 
to the corresponding cobordisms. The main difference to the Khovanov case
is that here not all faces are commutative. More precisely, the 
  square corresponding to the Torus relation $(T1)$ is anti--commutative. 
However, by the definition of the $\Z_2$--extension,
 the 2--cochain $\psi\in H^2(B,\mathbb{Z}/2\mathbb{Z})$ ($B$ is the hypercube) which associates $1$ to any anticommutative face of the 
hypercube and $-1$ to any commutative one is a cocycle, i.e.
it vanishes on the boundary of any cube. Since the
cube is contractible, any cocycle is a coboundary.
Consequently, there exists a function
$\ep:{\cal E}\to\Z_2$ from the set
of edges of the hypercube to $\Z_2$, called
 a sign assignment, such that $\ep(e_1)\ep(e_2)\ep(e_3)\ep(e_4)=\psi(D)$
 for any four edges $e_1,\dots,e_4$   forming a square $D$.
Hence,  multiplying edges of the hypercube by
the  signs  $\ep$, we get
 a chain complex
$(F\ll D\rr, \partial_\ep)$. It is easy to see
that this complex is independent on the choice of a sign
assignment.

\begin{lem}\label{22}
Given two sign assignments $\ep$ and $\ep'$, the chain complexes
$(F\ll D\rr, \partial_\ep)$
$(F\ll D\rr, \partial_{\ep'})$ are isomorphic.
\end{lem}

\begin{proof}
The product $\ep\ep'$ is a 1--cocycle. Since the hypercube is contractible,
 this 1-cocycle is a coboundary of
a 0--cochain $\eta: {\cal V}\to\Z_2$. The identity map times $\eta$
provides  the required isomorphism.
\end{proof}

 In the case, when $t=0$, 
the
homology groups of $(F\ll D\rr, \partial)$
 are graded and the graded Euler characteristic
coincides with the Jones polynomial. If $t\neq 0$, then
$\gr$ defines a filtration on our chain complex, 
similar to the one considered by Lee \cite{Lee}.

Our next aim is to show that the complex we just constructed  is
isomorphic to the Khovanov  complex.

\subsection{Proof of Theorem \ref{main}}
We have to show that  $(F\ll D\rr, \partial)$ and
$(F_{\Kh}\ll D\rr, \partial_{\Kh})$ are isomorphic.  

For any circle $c$ in $S^{\amalg k}$, we define $\deg(c)$ 
to be the number of circles in $S^{\amalg k}$  containing $c$ inside.
Further, we define an endomorphism $\phi_k$ of $ F(S^{\amalg k})$  as follows: 
For a copy of $A$ associated with $c$, we put
\be
\phi_c:\left\{\begin{array}{l@{\quad\mapsto\quad}r}
 \one & \one \\
 X&  (-1)^{\deg(c)} \; X\\
\end{array}
\right.
\ee
Then $\phi_k$ is the composition of $\phi_c$ for all circles in $S^{\amalg k}$.
By abuse of notation $\phi_k$ depends not only on $k$ but
also on the configuration of circles in $S^{\amalg k}$.

Given a link diagram $D$ with $d$ crossings, consider two 
Khovanov's hypercubes of resolutions associated with $D$.
Apply $F$ to one of them and $F_{\Kh}$ to 
the other and do not use any sign assignment, i.e.
all squares in $F_{\Kh} \ll D\rr$ are commutative.
Further, observe that
 with each vertex of the hypercube, there is a copy of $A^{\otimes k}$, for a certain $k$, associated.
Applying $\phi_{ k}$ to any such vertex, 
we get a  map $\Phi$ with the source
$(F\ll D\rr, \partial)$  and the target  $(F_{\Kh}\ll D\rr, \partial_{\Kh})$, without any sign assignments.
Our next goal is to see that there exists a sign assignment on the $(d+1)$--dimensional hypercube $\ll D\rr \times [0,1]$
 making $\Phi$ to a chain map.


For this, it is enough to check
 that each $3$--dimensional cube in this $(d+1)$--dimensional hypercube 
contains an even number of anticommutative faces. Note that there are three different cases:
 (1) the cube is contained in the source hypercube, (2) the cube is contained in the target hypercube,
 (3) the cube contains exactly one face in the source hypercube and one face in the target hypercube. 
The first case follows from Lemma \ref{cocy} and the second from the fact that $F_{\Kh}$ is a $(1+1)$ TQFT. 
The third case rely on a case by case check. Note that all  faces in 
the source hypercube correspond to  relations in $\NC_1$. Hence, we have to check the claim for any cube,
whose upper face is a relation in $\NC_1$, the lower face is the corresponding Khovanov square and whose  vertical
edges are labeled by $\Phi$. In addition, since the map $\phi_c$ depends on $\deg(c)$ explicitly, we
have to ensure that the claim holds after changing 
the nestedness of each circle by one. The tables below show
 that any cube of type (3) does 
 have only commutative or anticommutative faces. It is left to the 
reader to check that all cubes of type (3) do have an even number 
of anticommutative faces.  For this one has to consider  
all cubes where the upper face corresponds to one of the relations
 in Definition \ref{nc}. Moreover,  each cube should be checked twice 
for  different nestedness modulo 2.\\
\begin{center}\includegraphics{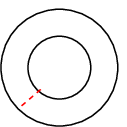} \hspace{0.5cm}
\begin{tabular}{|c|c|c|c|c|} \hline & $m_1$ & $\phi\circ m_1$ & $\phi$ & $m_0\circ \phi$\\
\hline $\one \otimes \one$ & $\one$ & $\one$ & $\one \otimes \one$ & $\one$\\
\hline$\one \otimes X$ & $X$ & $X$ & $\one \otimes X$ & $X$\\
\hline $X \otimes \one$ & $-X$ & $-X$ & $-X \otimes \one$ & $-X$\\
\hline $X \otimes X$ & $-t$ & $-t$ & $-X \otimes X$ & $-t$\\
\hline
\end{tabular}
\end{center}
\begin{center}\includegraphics{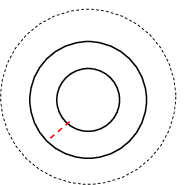} \hspace{0.5cm}
\begin{tabular}{|c|c|c|c|c|}\hline & $m_1$ & $\phi\circ m_1$ & $\phi$ & $m_0\circ \phi$\\
\hline $\one \otimes \one$ & $\one$ & $\one$ & $\one \otimes \one$ & $\one$\\
\hline $\one \otimes X$ & $X$ & $-X$ & $-\one \otimes X$ & $-X$\\
\hline $X \otimes \one$ & $-X$ & $X$ & $X \otimes \one$ & $X$\\
\hline $X \otimes X$ & $-t$ & $-t$ & $-X \otimes X$ & $-t$\\
\hline 
\end{tabular}
\end{center}
\begin{center} \includegraphics{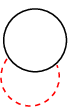} \hspace{0.5cm}
\begin{tabular}{|c|c|c|c|c|} \hline  & $\Delta_1$ & $\phi\circ \Delta_1$ & $\phi$ & $\Delta_0\circ \phi$\\
\hline $\one $ & $X \otimes \one-\one \otimes X$ & $-X\otimes \one -\one \otimes X$ & $\one$ & $X \otimes \one+\one \otimes X$\\
\hline $ X$ &$X\otimes X -t\one \otimes \one$  & $-X\otimes X -t \one \otimes \one $ & $X$ & $X\otimes X +t\one \otimes \one$\\
\hline 
\end{tabular}
\end{center}
\begin{center}\includegraphics{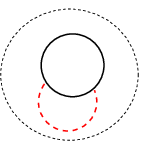} \hspace{0.5cm}
\begin{tabular}{|c|c|c|c|c|}\hline   & $\Delta_1$ & $\phi\circ \Delta_1$ & $\phi$ & $\Delta_0\circ \phi$\\
\hline $\one $ & $X \otimes \one-\one \otimes X$ & $X\otimes \one +\one \otimes X$ & $\one$ & $X \otimes \one+\one \otimes X$\\
\hline $ X$ &$X\otimes X -t\one \otimes \one$  & $-X\otimes X -t \one \otimes \one $ & $-X$ & $-X\otimes X -t\one \otimes \one$\\
\hline 
\end{tabular}
\end{center}
\vspace{0.5cm}
\begin{center} \includegraphics{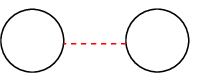} \hspace{0.5cm}
\begin{tabular}{|c|c|c|c|c|}\hline  & $m_0$ & $\phi\circ m_0$ & $\phi$ & $m_0\circ \phi$\\
\hline $\one \otimes \one$ & $\one$ & $\one$ & $\one \otimes \one$ & $\one$\\
\hline $\one \otimes X$ & $X$ & $X$ & $\one \otimes X$ & $X$\\
\hline $X \otimes \one$ & $X$ & $X$ & $X \otimes \one$ & $X$\\
\hline $X \otimes X$ & $t$ & $t$ & $X \otimes X$ & $t$\\
\hline 
\end{tabular}
\end{center}
\vspace{0.5cm}
\begin{center}\includegraphics{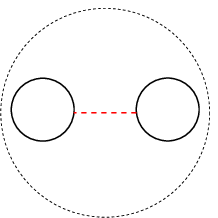} \hspace{0.5cm}
\begin{tabular}{|c|c|c|c|c|}\hline  & $m_0$ & $\phi\circ m_0$ & $\phi$ & $m_0\circ \phi$\\
\hline $\one \otimes \one$ & $\one$ & $\one$ & $\one \otimes \one$ & $\one$\\
\hline $\one \otimes X$ & $X$ & $-X$ & $-\one \otimes X$ & $-X$\\
\hline $X \otimes \one$ & $X$ & $-X$ & $-X \otimes \one$ & $-X$\\
\hline $X \otimes X$ & $+t$ & $+t$ & $+X \otimes X$ & $+t$\\
\hline 
\end{tabular}
\end{center}
\vspace{0.5cm}
\begin{center} \includegraphics{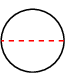} \hspace{0.5cm}
\begin{tabular}{|c|c|c|c|c|}\hline   & $\Delta_0$ & $\phi\circ \Delta_0$ & $\phi$ & $\Delta_0\circ \phi$\\
\hline $\one $ & $X \otimes \one+\one \otimes X$ & $X\otimes \one +\one \otimes X$ & $\one$ & $X \otimes \one+\one \otimes X$\\
\hline $ X$ &$X\otimes X +t\one \otimes \one$  & $X\otimes X +t \one \otimes \one $ & $X$ & $X\otimes X +t\one \otimes \one$\\
\hline 
\end{tabular}
\end{center}
\vspace{0.5cm}
\begin{center}\includegraphics{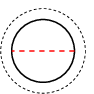} \hspace{0.5cm}
\begin{tabular}{|c|c|c|c|c|}\hline   & $\Delta_0$ & $\phi\circ \Delta_0$ & $\phi$ & $\Delta_0\circ \phi$\\
\hline $\one $ & $X \otimes \one+\one \otimes X$ & $-X\otimes \one -\one \otimes X$ & $\one$ & $X \otimes \one+\one \otimes X$\\
\hline $ X$ &$X\otimes X +t\one \otimes \one$  & $X\otimes X +t \one \otimes \one $ & $-X$ & $-X\otimes X -t\one \otimes \one$\\
\hline 
\end{tabular}
\end{center}
\vspace{0.5cm}

To finish, observe that the map $\phi$ composed with a sign assignment is clearly invertible, and hence, is the desired isomorphism.
\qed

\subsection{Proof of Theorem \ref{main1}}
Let us  search for further strict $\Z_2$--extensions of $\NC$ systematically.
For $e_i \in \Z_2$ with $i=1,...,10$, we put
\be
m_0:\left\{\begin{array}{l@{\quad\mapsto\quad}r}
\one\otimes \one & e_1\one\\
\one\otimes X& e_2 X\\
X\otimes \one& e_2 X\\
X\otimes X& 0
\end{array}
\right.
\quad\quad
m_1:\left\{\begin{array}{l@{\quad\mapsto\quad}r}
\one\otimes \one & e_5\one\\
\one\otimes X& e_6 X\\
X\otimes \one& e_7X\\
X\otimes X& 0
\end{array}
\right.
\ee

\be
\Delta_0:\left\{\begin{array}{l@{\quad\mapsto\quad}l}
 \one & e_3 (X\otimes\one+   \one\otimes X)\\
 X& 
e_4 X\otimes X \\
\end{array}
\right.
\quad\quad
\Delta_1:\left\{\begin{array}{l@{\quad\mapsto\quad}l}
 \one & e_8 X\otimes\one+  e_9 \one\otimes X\\
 X& 
e_{10} X\otimes X \\
\end{array}
\right.
\ee
The relations in Definition \ref{nc} should hold up to sign
for any EQFT. They impose the following relations on $e_i$:

4. row Frobenius type relations (1) $\Longrightarrow\;\;\;\;$ $e_6=e_5$,
$e_9=e_{10}$;

3. row Frobenius type relations (1) $\Longrightarrow\;\;\;\;$ $e_7 e_8=e_5 e_9$;

the ordinary Frobenius relation  $\Longrightarrow\;\;\;\;$
$e_1=e_2$, $e_3=e_4$.

Modulo these identities, 
there are 5 free parameters, i.e. 32 cases to consider.
It is a simple check that all of them produce the Khovanov
or nested Khovanov Frobenius system, after changing the sign of
one or two operations. 

It remains to construct
an isomorphism between, say, nested Khovanov complex 
and the one where  $m_0$ is replaced by $-m_0$. Let us consider the
map between two nested Khovanov hypercubes which is
identity on all vertices, 
except of the tails of edges corresponding to $m_0$,
at those edges the map is minus the identity. As in the proof of Theorem
\ref{main}, the cone of this map is a hypercube of a dimension one bigger.
Let convince our self that all 
3--dimensional cubes of that hypercube
have an even number of anticommutative faces.
We have to check only  cubes whose upper horizontal faces belong 
to the  nested
Khovanov complex, the bottom horizontal face to
 the nested Khovanov with $m_0$ replaced by $-m_0$,
  and whose vertical edges are given by our map.
If the upper horizontal face has an even number of $m_0$ maps, then the cube
has an even number of vertical anticommutative faces.
If it has an odd number of $m_0$ maps, then there is an odd number of
vertical anticommutative faces, but 
either the top or the
 bottom face is anticommutative.
Hence, like in the proof of Theorem \ref{main}, there exists a sign
arrangement on this hypercube providing the desired isomorphism.
\qed

\section{Odd Khovanov homology}


\subsection{The extension $\Odd$}
\begin{dfn} The extension $\Odd$ of $\Cob$ is
 defined as follows:
{\rm The objects of $\Odd$ are finite ordered
 set of circles. The morphisms are generated by
\begin{center}
\begin{picture}(0,0)%
\epsfig{file=gencob.pstex}%
\end{picture}%
\setlength{\unitlength}{789sp}%
\begingroup\makeatletter\ifx\SetFigFont\undefined%
\gdef\SetFigFont#1#2#3#4#5{%
  \reset@font\fontsize{#1}{#2pt}%
  \fontfamily{#3}\fontseries{#4}\fontshape{#5}%
  \selectfont}%
\fi\endgroup%
\begin{picture}(19012,4860)(7477,-3622)
\put(18001,-3361){\makebox(0,0)[lb]{\smash{{\SetFigFont{12}{14.4}{\familydefault}{\mddefault}{\updefault}$death$}}}}
\put(7651,-3361){\makebox(0,0)[lb]{\smash{{\SetFigFont{12}{14.4}{\familydefault}{\mddefault}{\updefault}$merge$}}}}
\put(11476,-3361){\makebox(0,0)[lb]{\smash{{\SetFigFont{12}{14.4}{\familydefault}{\mddefault}{\updefault}$birth$}}}}
\put(14851,-3361){\makebox(0,0)[lb]{\smash{{\SetFigFont{12}{14.4}{\familydefault}{\mddefault}{\updefault}$split$}}}}
\put(21151,-3361){\makebox(0,0)[lb]{\smash{{\SetFigFont{12}{14.4}{\familydefault}{\mddefault}{\updefault}$permutation$}}}}
\end{picture}%

\end{center}

subject to the 
following sets of relations:


(1) Commutativity and co--commutativity relation
\begin{center}
\begin{picture}(0,0)%
\epsfig{file=com.pstex}%
\end{picture}%
\setlength{\unitlength}{789sp}%
\begingroup\makeatletter\ifx\SetFigFont\undefined%
\gdef\SetFigFont#1#2#3#4#5{%
  \reset@font\fontsize{#1}{#2pt}%
  \fontfamily{#3}\fontseries{#4}\fontshape{#5}%
  \selectfont}%
\fi\endgroup%
\begin{picture}(11073,3273)(3427,-9010)
\put(9676,-7711){\makebox(0,0)[lb]{\smash{{\SetFigFont{12}{14.4}{\familydefault}{\mddefault}{\updefault}$=$}}}}
\end{picture}%
\hspace{2cm}
\begin{picture}(0,0)%
\epsfig{file=cocomminus.pstex}%
\end{picture}%
\setlength{\unitlength}{789sp}%
\begingroup\makeatletter\ifx\SetFigFont\undefined%
\gdef\SetFigFont#1#2#3#4#5{%
  \reset@font\fontsize{#1}{#2pt}%
  \fontfamily{#3}\fontseries{#4}\fontshape{#5}%
  \selectfont}%
\fi\endgroup%
\begin{picture}(11073,3273)(14602,-9010)
\put(19426,-7711){\makebox(0,0)[rb]{\smash{{\SetFigFont{12}{14.4}{\familydefault}{\mddefault}{\updefault}$=-$}}}}
\end{picture}%

\end{center}

(2) Associativity and coassociativity relations
\begin{center}
\begin{picture}(0,0)%
\epsfig{file=ass.pstex}%
\end{picture}%
\setlength{\unitlength}{789sp}%
\begingroup\makeatletter\ifx\SetFigFont\undefined%
\gdef\SetFigFont#1#2#3#4#5{%
  \reset@font\fontsize{#1}{#2pt}%
  \fontfamily{#3}\fontseries{#4}\fontshape{#5}%
  \selectfont}%
\fi\endgroup%
\begin{picture}(14748,5148)(18277,-10135)
\put(26326,-7861){\makebox(0,0)[rb]{\smash{{\SetFigFont{12}{14.4}{\familydefault}{\mddefault}{\updefault}$=$}}}}
\end{picture}%
\hspace{2cm}
\begin{picture}(0,0)%
\epsfig{file=coassminus.pstex}%
\end{picture}%
\setlength{\unitlength}{789sp}%
\begingroup\makeatletter\ifx\SetFigFont\undefined%
\gdef\SetFigFont#1#2#3#4#5{%
  \reset@font\fontsize{#1}{#2pt}%
  \fontfamily{#3}\fontseries{#4}\fontshape{#5}%
  \selectfont}%
\fi\endgroup%
\begin{picture}(14748,5148)(3427,-10135)
\put(10126,-7861){\makebox(0,0)[lb]{\smash{{\SetFigFont{12}{14.4}{\familydefault}{\mddefault}{\updefault}$=-$}}}}
\end{picture}%

\end{center}

(3) Frobenius relations
\begin{center}
\begin{picture}(0,0)%
\epsfig{file=frob.pstex}%
\end{picture}%
\setlength{\unitlength}{789sp}%
\begingroup\makeatletter\ifx\SetFigFont\undefined%
\gdef\SetFigFont#1#2#3#4#5{%
  \reset@font\fontsize{#1}{#2pt}%
  \fontfamily{#3}\fontseries{#4}\fontshape{#5}%
  \selectfont}%
\fi\endgroup%
\begin{picture}(21348,5165)(18277,-13068)
\put(25801,-10786){\makebox(0,0)[rb]{\smash{{\SetFigFont{12}{14.4}{\familydefault}{\mddefault}{\updefault}$=$}}}}
\put(33526,-10786){\makebox(0,0)[rb]{\smash{{\SetFigFont{12}{14.4}{\familydefault}{\mddefault}{\updefault}$=$}}}}
\end{picture}%

\end{center}

(4) Unit and Counit relations
\begin{center}
\begin{picture}(0,0)%
\epsfig{file=unitlaw.pstex}%
\end{picture}%
\setlength{\unitlength}{789sp}%
\begingroup\makeatletter\ifx\SetFigFont\undefined%
\gdef\SetFigFont#1#2#3#4#5{%
  \reset@font\fontsize{#1}{#2pt}%
  \fontfamily{#3}\fontseries{#4}\fontshape{#5}%
  \selectfont}%
\fi\endgroup%
\begin{picture}(11448,3222)(3577,-8569)
\put(9901,-7261){\makebox(0,0)[lb]{\smash{{\SetFigFont{12}{14.4}{\familydefault}{\mddefault}{\updefault}$=$}}}}
\end{picture}%
\hspace{2cm}
\begin{picture}(0,0)%
\epsfig{file=counitlaw.pstex}%
\end{picture}%
\setlength{\unitlength}{789sp}%
\begingroup\makeatletter\ifx\SetFigFont\undefined%
\gdef\SetFigFont#1#2#3#4#5{%
  \reset@font\fontsize{#1}{#2pt}%
  \fontfamily{#3}\fontseries{#4}\fontshape{#5}%
  \selectfont}%
\fi\endgroup%
\begin{picture}(11448,3222)(15127,-8569)
\put(20251,-7261){\makebox(0,0)[rb]{\smash{{\SetFigFont{12}{14.4}{\familydefault}{\mddefault}{\updefault}$=$}}}}
\end{picture}%

\end{center}

(5) Permutation relations
\begin{center}
\begin{picture}(0,0)%
\epsfig{file=trans2.pstex}%
\end{picture}%
\setlength{\unitlength}{789sp}%
\begingroup\makeatletter\ifx\SetFigFont\undefined%
\gdef\SetFigFont#1#2#3#4#5{%
  \reset@font\fontsize{#1}{#2pt}%
  \fontfamily{#3}\fontseries{#4}\fontshape{#5}%
  \selectfont}%
\fi\endgroup%
\begin{picture}(10698,3216)(7102,-8944)
\put(13126,-7561){\makebox(0,0)[lb]{\smash{{\SetFigFont{12}{14.4}{\familydefault}{\mddefault}{\updefault}$=$}}}}
\end{picture}%
\hspace{2cm}
\begin{picture}(0,0)%
\epsfig{file=R3trans.pstex}%
\end{picture}%
\setlength{\unitlength}{789sp}%
\begingroup\makeatletter\ifx\SetFigFont\undefined%
\gdef\SetFigFont#1#2#3#4#5{%
  \reset@font\fontsize{#1}{#2pt}%
  \fontfamily{#3}\fontseries{#4}\fontshape{#5}%
  \selectfont}%
\fi\endgroup%
\begin{picture}(19698,5166)(3952,-10294)
\put(12976,-8161){\makebox(0,0)[lb]{\smash{{\SetFigFont{12}{14.4}{\familydefault}{\mddefault}{\updefault}$=$}}}}
\end{picture}%

\end{center}

(6) Unit--Permutations and Counit--Permutation relations
\begin{center}
\begin{picture}(0,0)%
\epsfig{file=permunit.pstex}%
\end{picture}%
\setlength{\unitlength}{789sp}%
\begingroup\makeatletter\ifx\SetFigFont\undefined%
\gdef\SetFigFont#1#2#3#4#5{%
  \reset@font\fontsize{#1}{#2pt}%
  \fontfamily{#3}\fontseries{#4}\fontshape{#5}%
  \selectfont}%
\fi\endgroup%
\begin{picture}(11073,3222)(4027,-7975)
\put(10201,-6661){\makebox(0,0)[lb]{\smash{{\SetFigFont{12}{14.4}{\familydefault}{\mddefault}{\updefault}$=$}}}}
\end{picture}%
\hspace{2cm}
\begin{picture}(0,0)%
\epsfig{file=permcounit.pstex}%
\end{picture}%
\setlength{\unitlength}{789sp}%
\begingroup\makeatletter\ifx\SetFigFont\undefined%
\gdef\SetFigFont#1#2#3#4#5{%
  \reset@font\fontsize{#1}{#2pt}%
  \fontfamily{#3}\fontseries{#4}\fontshape{#5}%
  \selectfont}%
\fi\endgroup%
\begin{picture}(11523,3222)(6952,-8200)
\put(14626,-6736){\makebox(0,0)[rb]{\smash{{\SetFigFont{12}{14.4}{\familydefault}{\mddefault}{\updefault}$=$}}}}
\end{picture}%

\end{center}

(7) Merge--Permutation and Split--Permutation relations
\begin{center}
\begin{picture}(0,0)%
\epsfig{file=permerge.pstex}%
\end{picture}%
\setlength{\unitlength}{789sp}%
\begingroup\makeatletter\ifx\SetFigFont\undefined%
\gdef\SetFigFont#1#2#3#4#5{%
  \reset@font\fontsize{#1}{#2pt}%
  \fontfamily{#3}\fontseries{#4}\fontshape{#5}%
  \selectfont}%
\fi\endgroup%
\begin{picture}(16698,5166)(18277,-10144)
\put(26101,-7861){\makebox(0,0)[rb]{\smash{{\SetFigFont{12}{14.4}{\familydefault}{\mddefault}{\updefault}$=$}}}}
\end{picture}%
\hspace{0.5cm}
\begin{picture}(0,0)%
\epsfig{file=permsplit.pstex}%
\end{picture}%
\setlength{\unitlength}{789sp}%
\begingroup\makeatletter\ifx\SetFigFont\undefined%
\gdef\SetFigFont#1#2#3#4#5{%
  \reset@font\fontsize{#1}{#2pt}%
  \fontfamily{#3}\fontseries{#4}\fontshape{#5}%
  \selectfont}%
\fi\endgroup%
\begin{picture}(16398,5166)(-6548,-10144)
\put(-449,-7861){\rotatebox{360.0}{\makebox(0,0)[lb]{\smash{{\SetFigFont{12}{14.4}{\familydefault}{\mddefault}{\updefault}$=$}}}}}
\end{picture}%

\end{center}
(8) Commutation relations
\begin{center}
\begin{picture}(0,0)%
\epsfig{file=splitsplitcomutation.pstex}%
\end{picture}%
\setlength{\unitlength}{789sp}%
\begingroup\makeatletter\ifx\SetFigFont\undefined%
\gdef\SetFigFont#1#2#3#4#5{%
  \reset@font\fontsize{#1}{#2pt}%
  \fontfamily{#3}\fontseries{#4}\fontshape{#5}%
  \selectfont}%
\fi\endgroup%
\begin{picture}(13848,6732)(3727,-11269)
\put(9901,-8086){\makebox(0,0)[lb]{\smash{{\SetFigFont{12}{14.4}{\familydefault}{\mddefault}{\updefault}$=-$}}}}
\end{picture}%
\hspace{2cm}
\begin{picture}(0,0)%
\epsfig{file=counitsplitcomutation.pstex}%
\end{picture}%
\setlength{\unitlength}{789sp}%
\begingroup\makeatletter\ifx\SetFigFont\undefined%
\gdef\SetFigFont#1#2#3#4#5{%
  \reset@font\fontsize{#1}{#2pt}%
  \fontfamily{#3}\fontseries{#4}\fontshape{#5}%
  \selectfont}%
\fi\endgroup%
\begin{picture}(13848,4947)(3727,-11269)
\put(9676,-9361){\makebox(0,0)[lb]{\smash{{\SetFigFont{12}{14.4}{\familydefault}{\mddefault}{\updefault}$=-$}}}}
\end{picture}%
\\
\end{center}
\vspace{1cm}
\begin{center}
\begin{picture}(0,0)%
\epsfig{file=counitcounitcomutation.pstex}%
\end{picture}%
\setlength{\unitlength}{789sp}%
\begingroup\makeatletter\ifx\SetFigFont\undefined%
\gdef\SetFigFont#1#2#3#4#5{%
  \reset@font\fontsize{#1}{#2pt}%
  \fontfamily{#3}\fontseries{#4}\fontshape{#5}%
  \selectfont}%
\fi\endgroup%
\begin{picture}(11748,3978)(3727,-10300)
\put(8776,-8611){\makebox(0,0)[lb]{\smash{{\SetFigFont{12}{14.4}{\familydefault}{\mddefault}{\updefault}$=-$}}}}
\end{picture}%

\end{center}
All the other commutation relations hold with plus sign.}\\
\label{odddef}
\end{dfn}
All axioms of a semistrict monoidal 2--category are satisfied.
\begin{remark}
For another definition of the
 semistrict monoidal just described, 
we endow the morphisms with the following $\mathbb{Z}/2\mathbb{Z}$--grading:
\be 
\rm{sd}\left(\begin{array}{c}\begin{picture}(0,0)%
\epsfig{file=merge.pstex}%
\end{picture}%
\setlength{\unitlength}{395sp}%
\begingroup\makeatletter\ifx\SetFigFont\undefined%
\gdef\SetFigFont#1#2#3#4#5{%
  \reset@font\fontsize{#1}{#2pt}%
  \fontfamily{#3}\fontseries{#4}\fontshape{#5}%
  \selectfont}%
\fi\endgroup%
\begin{picture}(3198,3198)(7477,-2110)
\end{picture}%
\end{array}\right)=0, \quad \rm{sd}\left(\begin{array}{c}\begin{picture}(0,0)%
\epsfig{file=split.pstex}%
\end{picture}%
\setlength{\unitlength}{395sp}%
\begingroup\makeatletter\ifx\SetFigFont\undefined%
\gdef\SetFigFont#1#2#3#4#5{%
  \reset@font\fontsize{#1}{#2pt}%
  \fontfamily{#3}\fontseries{#4}\fontshape{#5}%
  \selectfont}%
\fi\endgroup%
\begin{picture}(3198,3198)(14227,-2035)
\end{picture}%
\end{array}\right)=1,
\ee
\be 
\rm{sd}\left(\begin{array}{c}\begin{picture}(0,0)%
\epsfig{file=birth.pstex}%
\end{picture}%
\setlength{\unitlength}{395sp}%
\begingroup\makeatletter\ifx\SetFigFont\undefined%
\gdef\SetFigFont#1#2#3#4#5{%
  \reset@font\fontsize{#1}{#2pt}%
  \fontfamily{#3}\fontseries{#4}\fontshape{#5}%
  \selectfont}%
\fi\endgroup%
\begin{picture}(1098,1278)(11977,-1075)
\end{picture}%
\end{array}\right)=0, \quad \rm{sd}\left(\begin{array}{c}\begin{picture}(0,0)%
\epsfig{file=death.pstex}%
\end{picture}%
\setlength{\unitlength}{395sp}%
\begingroup\makeatletter\ifx\SetFigFont\undefined%
\gdef\SetFigFont#1#2#3#4#5{%
  \reset@font\fontsize{#1}{#2pt}%
  \fontfamily{#3}\fontseries{#4}\fontshape{#5}%
  \selectfont}%
\fi\endgroup%
\begin{picture}(1098,1278)(18802,-1075)
\end{picture}%
\end{array}\right)=1.
\ee
 This grading is additive under composition and disjoint union.
 The  monoidal structure $\boxtimes$ on $\Odd$  
 can be  defined as follows:\\
For any two generators $f$ and $g$ and the permutation $\rm{Perm}$, 
\be f\boxtimes \id:=f\otimes \id,\quad
f\boxtimes g:=(f\boxtimes \id)\circ(\id \boxtimes g) \ee
where $\otimes$ denotes the disjoint union.
The composition rule is modified as follows:
\be
(\id \boxtimes g)\circ(f \boxtimes \id)=
{(-1)}^{\rm{sd}(f)\;\rm{sd}(g)}\; f\boxtimes g.
\ee
\end{remark}

For an alternative description of $\Odd$, see Putyra's Master Thesis 
\cite{Pu} using cobordisms with chronology.\\

Let us  check that $\Odd$ is indeed an extension.
\begin{lem}
The automorphism group of any 1--morphism is trivial.
\end{lem}
\begin{proof}
The relations imply that all squares depicted on the RHS, resp. LHS, of
Figure 2 in \cite{Odd}
 are commutative (resp. anticommutative).
The result follows now from Lemma 2.1 in \cite{Odd}, showing
 that any cube has even number of anticommutative faces and additional checks like the one in the relations satisfied by the 2--morphisms in Lemma 32 \cite{Lauda}.
The result can also be checked completely by hand by proving that all the relations in Lemma 32 \cite{Lauda} are satisfied by the 2--morphisms in $\Odd$ which are only signs.
Many of them are obvious, since many 2--morphisms in our case are just identities. The fact that this is enough still follows from Bergman's Diamond lemma \cite{Berg}. 

\end{proof}

\subsection{Odd Frobenius system}
In \cite{Odd}, an EQFT into the $\mathbb{Z}/2\mathbb{Z}$--graded abelian groups based on $\Odd$
 is constructed. 

Using Khovanov's algebra $A_0=\Z[X]/X^2$, one can describe this 
 EQFT $F_0: \Odd\to \Z$-$\mod$ as follows: 
$F_0$ maps a circle to $A$ where $A$ is $\mathbb{Z}/2\mathbb{Z}$--graded as follows:  
$\one$ is in degree $0$ and $X$ is in degree $1$. To $n$ circles, $F_0$ assigns
 $A^{\otimes n}$.
To generating morphisms $F_0$ assigns the following maps:
\be
m:\left\{\begin{array}{l@{\quad\mapsto\quad}r}
\one\otimes \one & \one\\
\one\otimes X& X\\
X\otimes \one&  X\\
X\otimes X& 0
\end{array}
\right.
\quad\quad\quad
P:\left\{\begin{array}{l@{\quad\mapsto\quad}r}
\one\otimes \one & \one\otimes \one\\
\one\otimes X& X\otimes \one\\
X\otimes \one& \one\otimes X\\
X\otimes X& -X\otimes X
\end{array}
\right.
\ee
$$
\Delta:\left\{\begin{array}{l@{\quad\mapsto\quad}l}
 \one & X\otimes\one  -   \one\otimes X\\
 X&  X\otimes X \\
\end{array}
\right.
$$
The  maps $\e$ and $\eta$ are the same as in Khovanov case.

Due to the fact that $\Delta$ and $\e$ are of degree $1$, $F_0$ can not map disjoint union of cobordisms to the tensor product
of  maps assigned to them, since in this case relations (6) and (7) would not be satisfied.
Instead, $F_0$ maps disjoint union to $\boxtimes$ defined
as follows:
\be f\boxtimes \id:=f\otimes \id,\quad
\id\boxtimes f:=\rm{Perm}\circ(f\otimes\id)\circ \rm{Perm},\quad
f\boxtimes g:=(f\boxtimes \id)\circ(\id \boxtimes g) \ee
\be
(\id \boxtimes g)\circ(f \boxtimes \id)={(-1)}^{\rm{deg}(f)\rm{deg}(g)}f\boxtimes g
\ee
The relations (6) and (7) hold now just by definition.

Applied to the Khovanov hypercube, this EQFT gives rise
to a  link homology theory, called odd Khovanov homology \cite{Odd}.

\end{document}